\documentstyle[12pt]{article}
\newcommand{\be}{\begin{equation}}
\newcommand{\ee}{\end{equation}}
\newcommand{\ba}{\begin{eqnarray}}
\newcommand{\ea}{\end{eqnarray}}

\newcommand{\lab}[1]{\label{#1}}
\newcommand{\re}[1]{(\ref{#1})}
\newcommand{\bb}{\begin{thebibliography}}
\newcommand{\eb}{\end{thebibliography}}

\newtheorem{tm}{Theorem}
\newtheorem{lem}{Lemma}

\begin{document}

\begin{center}{{\LARGE {Biorthogonality of the Lagrange interpolants}}}

\vspace{5mm}

{\large \bf Alexei Zhedanov}\\

\vspace{2mm}

{\it Donetsk Institute for Physics and Technology, Donetsk 83114, Ukraine}

\end{center}


\begin{abstract}
We show that the Lagrange interpolation polynomials are biorthogonal with
respect to a set of rational functions whose poles coinicde with interpolation points

\end{abstract}

The Newton-Lagrange interpolation is a well-known problem
in elementary calculus. Recall basic facts concerning this problem \cite{MT}, \cite{Gel}.

Let $A_k, \; k=0,1, 2,\dots$ and $a_k, \; k=0,1,2 ,\dots$  be two arbitrary sequences
of complex numbers (we assume that all $a_k$ are distinct $a_k \ne a_j$ if $ k \ne j$.
By interpolation polynomial we mean a $n$-degree polynomial $P_n(z)$
whose values at points $a_0, a_1, \dots, a_{n}$ coincide with $A_0, A_1, \dots, A_{n}$,
i.e.
\be
P_n(a_k) = A_k, \quad k=0, 1,2,\dots, n
\lab{i_cond} \ee
Usually the parameters $A_k$ are interpreted as values of some function $F(z)$
at fixed points $a_k$, i.e.
\be
A_k = F(a_k)
\lab{AFk} \ee
In this case polynomials $P_n(z)$ interpolate the function $F(z)$ at points $a_k$.

Explicit expression for interpolation polynomial $P_n(z)$ can be presented
in two forms. In the Newtonian form we have \cite{MT}, \cite{Gel}
\be
P_n(z) = \sum_{k=0}^{n} [a_0,a_1, \dots, a_k  ] \omega_k(x),
\lab{N_P} \ee
where
$$\omega_0=1, \; \omega_k(x) = (x-a_0)(x-a_1) \dots (x-a_{k-1})$$
and $[a_0,a_1, \dots, a_k  ]$ denotes the $k$-th Newtonian divided
differnce which is defined as
$$
[a_0]=A_0, \; [a_0,a_1]= \frac{A_1-A_0}{a_1-a_0}, \dots, [a_0, a_1, \dots, a_k]= \sum_{s=0}^k \frac{A_s}{\omega'_{k+1}(a_s)},
$$
where
$$
\omega_{k+1}'(a_s) = (a_s -a_1) (a_s-a_2) \dots (a_s - a_{s-1})(a_s - a_{s+1}) \dots (a_s - a_k) =
\prod_{i=0, i \ne s }^k(a_s-a_i)
$$

If representation \re{AFk} holds where $F(z)$ is a meromorphic function then
the Hermite formula is useful \cite{Gel}
\be
[a_0,a_1, \dots, a_k  ] =(2 \pi i)^{-1} \int_{\Gamma}
\frac{F(\zeta) d \zeta}{\omega_{k+1}(\zeta)},
\lab{D_Her} \ee
where the contour $\Gamma$ in complex plane is chosen such that points
$a_0, a_1, \dots, a_k$ lie inside the contour whereas all singularities
of the function $F(z)$ lie outside the contour.

In the Lagrangean form we have \cite{MT}, \cite{Gel}
\be
P_n(z) = \omega_{n+1}(z) \: \sum_{k=0}^n \frac{A_k}{(z-a_k) \omega_{n+1}'(a_k)},
\lab{L_P} \ee

Introduce the monic interpolation polynomials $\hat P_n(z) = P_n(z)/\alpha_n$, where
$$
\alpha_n =[a_0,a_1, \dots, a_n]
$$
In what follows we will assume that $\alpha_n \ne 0$ for all $n=1,2, \dots$.
This condition
guarantees that polynomials $P_n(z)$ are indeed of the $n$-th degree.
It is easily seen that $P_n(z) = z^n + O(z^{n-1})$.
For polynomials $\hat P_n(z)$
one has the recurrence relation \cite{Henr}
\be
\hat P_{n+1}(z) = (z-a_n + \frac{\alpha_n}{\alpha_{n+1}}) \hat P_n(z) -
\frac{\alpha_{n-1}}{\alpha_{n}} (z-a_n) \hat P_{n-1}(z)
\lab{rec_P} \ee
with the initial conditions
\be
\hat P_{-1}=0, \; \hat P_0(z) =1
\lab{ini_P} \ee
It is clear that the set of monic interpolation polynomials
$\hat P_0, \hat P_1(z), \dots, \hat P_n(z), \dots$
does not belong to a set of orthogonal polynomials (OP), becuase OP satisfy
3-term recurrence realtions of the form \cite{Sz}
\be
P_{n+1}(z) + b_n P_n(z) + u_n P_{n-1}(z) = z P_n(z)
\lab{OP} \ee
which doesn't have the form \re{rec_P}.

Nevertheless, recurrence relation \re{rec_P}
belongs to the so-called $R_I$-type recurrence relations (in terminology
of \cite{IM}). It was shown in \cite{IM} that polynomials satisfying
$R_I$-type relations possess some orthogonality property. In our case
this orthogonality property is well known \cite{MT}:

\begin{lem}
Polynomials $\hat P_n(z)$ satisfy formal orthogonality relation
\be I_{nj}=(2 \pi i)^{-1}\int_{\Gamma} \frac{\zeta^j \hat
P_n(\zeta) d \zeta}{\omega_{n+1}(\zeta) F(\zeta)} =
\frac{\delta_{nj}}{\alpha_n}, \quad j=0,1,\dots, n \lab{ort_G1}
\ee where the contour $\Gamma$ encompasses points $a_0, a_1,
\dots, a_N$ with $N \ge n$ and all singularities of the function
$1/F(z)$ lie outside the contour.
\end{lem}
For the proof it is sufficient to note, that under conditions upon choice of
the contour $\Gamma$, the intergral can be presented as a sum of residues
\be
I_{nj}= \sum_{s=0}^n \frac{a_s^j \hat P_n(a_s)}{A_s \omega'_{n+1}(a_s)} =
\sum_{s=0}^n \frac{a_s^j }{\alpha_n \omega'_{n+1}(a_s)},
\lab{Inj} \ee
where we used interpolaiton property \re{i_cond}. Hence, in the integral \re{ort_G1}
one can replace $\frac{P_n(\zeta)}{F(\zeta)}=1/\alpha_n$ and we have
$$
I_{nj}= (2 \pi i)^{-1}\int_{\Gamma} \frac{\zeta^j d \zeta}{\alpha_n \omega_{n+1}(\zeta)} =0, \; j=0,1,\dots, n-1
$$
because the value of the integral from pure rational function doesnot depend on choice of the contour $\Gamma$
(provided that all poles of the function lie inside the contour) and we can
choose $\Gamma$ as a circle of a great radius. For $j=n$ we have analogously
$$
I_{nn} = (2 \pi i)^{-1}\int_{\Gamma} \frac{\zeta^n d \zeta}{\alpha_n \omega_{n+1}(\zeta)} =
(2 \pi i)^{-1}\int_{\Gamma} \frac{d \zeta}{\alpha_n \zeta} = 1/\alpha_n.
$$
As was shown in \cite{ZheR} , \cite{LZ} polynomials of $R_I$ type possesses not only orthogonality
of the form \re{ort_G1} but also nice {\it biorthogonality} property with respect to
some set of rational functions.

In order to derive this biorthogonality propety in our case, we construct
auxiliary polynomials
\be
T_n(z) = \hat P_{n+1}(z) - (z-a_{n+1}) \hat P_n(z).
\lab{Tn} \ee

Clearly,
degree of these polynomials $\le n$. More exactly,
$$
T_n(z) = \nu_n z^n + O(z^{n-1}),
$$
where
\be
\nu_n=a_{n+1} -a_n + \frac{\alpha_n}{\alpha_{n+1}} -\frac{\alpha_{n-1}}{\alpha_{n}}
\lab{lead_T} \ee

In what folows we will assume that $\nu_n \ne 0$. This means that degree of $T_n(z)$ is $n$ and
it is possible to introduce monic polynomials
\be
\hat T_n = T_n(z)/\nu_n = z^n + O(z^{n-1})
\lab{monT} \ee

We have
\begin{tm}
Let $\hat P_n(z)$ be Lagrange interpolation polynomials for the
function $F(z)$ and $\hat T_n(z)$ defined by \re{Tn}, \re{monT}.
Define a set of rational functions \be V_n(z) = \frac{\hat
T_n(z)}{\omega_{n+2}(x)} \lab{Vn} \ee Then Lagrange interpolation
polynomials $\hat P_n(z)$ and functions $V_n(z)$ form a
biortogonal system in the folowing sense: \be \int_{\Gamma}
\frac{\hat P_n(\zeta) V_m(\zeta) d \zeta}{F(\zeta)} = \alpha_n
^{-1} \delta_{nm}, \lab{biortPV} \ee where the contour $\Gamma$
should be chosen such that  interpolation points $a_0, a_1, \dots,
a_N$ lie inside the contour $(N \ge \max(n,m+1))$, and the
function $1/F(z)$ is regular inside and on the contour.
\end{tm}
{\it Proof}. Assume first that $m<n-1$. Then we have, obviously
\ba
&&\int_{\Gamma} \frac{\hat P_n(\zeta) V_m(\zeta) d \zeta}{F(\zeta)} =
\int_{\Gamma} \frac{\hat P_n(\zeta) \hat T_m(\zeta)(z-a_{m+2})(z-a_{m+3}) \dots (z- a_{n})  d \zeta}{\omega_{n+1}(z) F(\zeta)}=
\nonumber \\
&&\int_{\Gamma} \frac{\hat P_n(\zeta) q_{n-1}(z) d \zeta }{F(\zeta)} =0
\lab{int1} \ea
where $q_{n-1}(z)$ is a polynomial of degree $\le n-1$ and in the last equality in
\re{int1} we used orthogonality property \re{ort_G1}.

If $m=n-1$ then
$$
\int_{\Gamma} \frac{\hat P_n(\zeta) V_{n-1}(\zeta) d \zeta}{F(\zeta)}=
\int_{\Gamma} \frac{\hat P_n(\zeta) \hat T_{n-1} (\zeta) d \zeta}{\omega_{n+1}(z)F(\zeta)}=0
$$
again by \re{ort_G1}

If $m>n$ then we can write down
$$
\int_{\Gamma} \frac{\hat P_n(\zeta) V_{m}(\zeta) d \zeta}{F(\zeta)}=
\int_{\Gamma} \frac{\hat P_n(\zeta) (\hat P_{m+1}(\zeta) -(\zeta-a_{m+1})\hat P_m(\zeta))  d\zeta}
{\omega_{m+2}(\zeta)F(\zeta)} =
$$
$$
\int_{\Gamma} \frac{\hat P_n(\zeta) \hat P_{m+1}(\zeta)  d\zeta}
{\omega_{m+2}(\zeta)F(\zeta)} - \int_{\Gamma} \frac{\hat P_n(\zeta) \hat P_m(\zeta)  d\zeta}
{\omega_{m+1}(\zeta)F(\zeta)} =0
$$
because both terms in the last relation vanish due to \re{ort_G1} for $n<m$.

Finaly, consider the case $m=n$. We have
$$
\int_{\Gamma} \frac{\hat P_n(\zeta) V_{n}(\zeta) d \zeta}{F(\zeta)}=
\int_{\Gamma} \frac{\hat P_n(\zeta) \hat T_{n} (\zeta) d \zeta}{\omega_{n+1}(\zeta)F(\zeta)}=
$$
$$
\int_{\Gamma} \frac{\hat P_n(\zeta) \hat P_{n+1} (\zeta) d \zeta}{\omega_{n+2}(\zeta)F(\zeta)} -
\int_{\Gamma} \frac{\hat P_n^2(\zeta) d \zeta}{\omega_{n+1}(\zeta)F(\zeta)} =
$$
$$
\int_{\Gamma} \frac{\hat P_n^2(\zeta) d \zeta}{\alpha_n \omega_{n}(z)} =1/\alpha_n
$$
The theorem is proven.

The biorthogonality property can be rewritten in another form \be
\sum_{s=0}^N  \lim_{z=a_s}{((z-a_s) P_n(z) V_m(z))}/A_s =
\delta_{nm}/\alpha_n, \lab{biort2} \ee where $N$ is any positive
integer such that $N \ge \max(n,m+1)$.

This result allows one to find coefficients $\xi_k$ in expansion of a given polynomial $Q(z)$ of degree $n$
in terms of interpolation polynomials $P_n(z)$
$$
Q(z) = \sum_{k=0}^n \xi_k \hat P_k(z)
$$
Indeed, from \re{biort2} we have
\be
\xi_k = (2 \pi i)^{-1} \alpha_k \int_{\Gamma} Q(\zeta) V_k(\zeta) d \zeta
\lab{exp_coef} \ee

Consider an example. For the exponential function $F(z) = \exp(h z)$
(with an arbitrary nonzero real parameter $h$) choose uniform grid of the
interpolation points $a_k =k, \; k=0,1,\dots$.
We then have (cf. \cite{Henr})
\be
P_n(z) = \sum_{k=0}^n \frac{(-z)_k}{k!}(1-e^h)^k
\lab{P_exp} \ee
where $(b)_k = b(b+1) \dots (b+k-1)$ is the Pochhammer symbol. From \re{P_exp} it is found
\be
\alpha_n = \frac{(e^h-1)^n}{n!}
\lab{alp_exp} \ee
Construct auxiliary polynomials
$T_n(z) = \hat P_{n+1}(z) - (z-a_{n+1}) \hat P_n(z)$. For the leading coefficient $\nu_n$
of the polynomials $T_n(z) = \nu_n z^n + O(z^{n-1})$
we have from \re{lead_T} and \re{alp_exp}
$$
\nu_n = \frac{e^h}{e^h-1} \ne 0
$$
So polynomials $T_n(z)$ are indeed of degree $n$ and for monic polynomials $\hat T_n (z) = T_n(z)/\nu_n$
it is not difficult to obtain a rather attractive closed foormula
\be
T_n(z) = \frac{(n+1)!}{(e^h-1)^n} \: {_2}F_1\left({-n, -z \atop -1-n}; 1-e^h \right)
\lab{T_exp} \ee
Thus for rational corresponding rational functions $V_n(z)$ we have from \re{Vn}
\be
V_n(z) = \frac{(n+1)!}{(1-e^h)^n (-z)_{n+2}} \:
{_2}F_1\left({-n, -z \atop -1-n}; 1-e^h \right)
\lab{V_exp} \ee
Using standard transformation formulas for the Gauss hypergeometric
function \cite{Bat}, we can present the functions $V_n(z)$ in a slightly different form
\be
V_n(z) = \frac{1}{(1-e^h)^n z(z-1)} \:
{_2}F_1\left({-n, -z \atop 2-z}; e^h \right)
\lab{V_exp2} \ee
Thus rational functions $V_n(z)$ form a biorthogonal set with
respect to the Lagrangean interpolaiton polynomials \re{P_exp}:
\be
\int_{\Gamma} P_n(\zeta) V_m(\zeta) \exp(-h \zeta) d \zeta = \delta_{nm}
\lab{bi_exp} \ee
where $\Gamma$ is an arbitrary contour containing the
points $0,1,\dots, \max(n,m+1)$ inside.

Note finally, that recurrence relation \re{rec_P} {\it completely
characterizes} the Lagrange interpolation polynomials $P_n(z)$. More exactly,
we have the
\begin{tm}
Assume that a set of monic $n$-th degree polynomials $\hat P_n(z)$
satisfies recurrence relation \re{rec_P} with initial conditions
\re{ini_P}, where parameters $\alpha_n, a_n, \; n=0,1,\dots$ are
arbitrary with the restrictions that all $a_i$ are distinct: $a_i
\ne a_j$, for $i \ne j$ and all $\alpha_n$ are nonzero $\alpha_n
\ne 0, \: n=0,1,\dots$. Then polynomials $P_n(z) = \alpha_n \hat
P_n(z)$ satisfy interpolation condition $P_n(a_k) = A_k, \;
k=0,1,\dots,n$ for all $n=0,1,\dots$, where \be A_n = P_n(a_n)=
\sum_{s=0}^n \alpha_s \omega_{s}(a_n), \; n=0, 1,2,\dots \lab{As}
\ee
\end{tm}
{\it Proof}. From recurrence relation \re{rec_P} and initial conditions \re{ini_P}
it can be easily found
\be
\hat P_{n+1}(z) - \frac{\alpha_n}{\alpha_{n+1}} \hat P_n(z) = \omega_{n+1}(z), \; n=0,1, \dots
\lab{rel1} \ee
Hence for $P_n(z) = \alpha_n \hat P_n(z)$ we have the conditions
$$
P_{n+1}(a_k) = P_n(a_k), \; k=0,1,\dots, n, \; n=0,1,2,\dots
$$
From these relations, by induction, we obtain
$$
P_n(a_k) = P_k(a_k)=A_k, \; k=0,1,\dots,n
$$
Thus interpolation conditions are fulfilled.
Expression \re{As} for $A_k$ follows then
from Newton formula \re{N_P}.

\end{document}